\documentclass[12pt,twoside]{amsart}
\usepackage{amsthm,amsfonts,amssymb,amscd,euscript}

\usepackage[matrix,arrow]{xy}

\author{Florin Ambro} 
\address{RIMS, Kyoto University\\
Kyoto 606-8502, Japan.}
\email{ambro@kurims.kyoto-u.ac.jp}

\newcommand{\Q}{{\mathbb Q}}
\newcommand{\Z}{{\mathbb Z}}

\newcommand{\R}{{\mathbb R}}

\newcommand{\calO}{{\mathcal O}}

\newcommand{\bA}{{\mathbf A}}
\newcommand{\bB}{{\mathbf B}}

\newcommand{\bD}{{\mathbf D}}
\newcommand{\bE}{{\mathbf E}}

\newcommand{\bH}{{\mathbf H}}
\newcommand{\bK}{{\mathbf K}}

\newcommand{\bM}{{\mathbf M}}

\newcommand{\bP}{{\mathbf P}}

\newcommand{\Alb}{\operatorname{Alb}}

\newcommand{\Exc}{\operatorname{Exc}}

\newcommand{\mult}{\operatorname{mult}}

\newcommand{\rank}{\operatorname{rank}}

\newcommand{\Spec}{\operatorname{Spec}}

\theoremstyle{plain}
\newtheorem{thm}{Theorem}[section]

\newtheorem{lem}[thm]{Lemma}

\newtheorem{prop}[thm]{Proposition}

\newtheorem{conj}[thm]{Conjecture}

\theoremstyle{definition}
\newtheorem{defn}[thm]{Definition}

\newtheorem{example}{Example}

\newtheorem{rem}[thm]{Remark}

\theoremstyle{remark}

\begin{document}

\bibliographystyle{amsalpha+}
\title[The nef dimension of log minimal models]
{The nef dimension of log minimal models}

\begin{abstract} This is the resume of the talk delivered 
by the author at the Symposium on Hodge theory, Degeneration 
and Complex surfaces, Tagajo, Miyagi, March 2004.
It is based on the papers~\cite{abnef,modbdiv}, with
some improvements.
\end{abstract}

\maketitle


\setcounter{section}{-1}


\section{Introduction}


\footnotetext[1]{This work was supported through a 
Twenty-First Century COE Kyoto Mathematics Fellowship.
}

The aim of this note is to discuss the Log 
Abundance Conjecture for log minimal models 
in terms of the nef dimension of their log 
canonical divisor. 

A {\em log minimal model} 
is a complex projective log variety $(X,B)$ 
with Kawamata log terminal singularities, such 
that the log canonical divisor $K+B$ is nef. 
The Log Abundance Conjecture~\cite{reid,KMM}
predicts that the pluricanonical linear system
$|m(K+B)|$ of a log minimal model is base point
free for some positive integer $m$. This conjecture
is closely related to the Log Minimal Model Program
Conjecture. A special case of Log Abundance, the 
Base Point Free Theorem~\cite{KMM}, played a key role
in establishing the known steps of the Log Minimal
Model Program (see~\cite{KMM}). 
Conversely, the $3$-dimensional case
of the Log Minimal Model Program was one of the key 
ingredients in the proof of Log Abundance in dimension
$3$ (\cite{ab,Miy1,Miy2,Miy3,Ka92,LAB,LAB'}).

The traditional approach towards the abundance problem 
is due to Kawamata~\cite{ab}. Modulo the Log Minimal
Model Program and Log Abundance in dimension at most
$d-1$, Log Abundance in dimension $d$ is equivalent 
to the following {\em non-vanishing statement}: 
given a $d$-dimensional log minimal model $(X,B)$, 
if there exists a curve $C$ in $X$ such that 
$(K+B)\cdot C>0$ then $\dim H^0(X,m(K+B))\ge 2$ for 
some positive integer $m$. 
This non-vanishing statement was established by Miyaoka 
and Kawamata in the case of minimal $3$-folds, and it was 
extended to the logarithmic case by Keel, Matsuki 
and McKernan. Their approach is based on a delicate analysis 
of the coefficients appearing in the Riemann-Roch formula 
for the pluricanonical divisors, and it is unclear whether
a similar approach could work in dimension $4$ or
higher.

In the following, we present a numerical approach towards 
the abundance problem, based on the nef dimension of a nef 
divisor, an invariant recently introduced by Tsuji~\cite{T00} 
and Bauer et al~\cite{nr}.
Given a nef divisor $D$ on a normal projective variety $X$, 
there exists a rational dominant map 
$f\colon X\dashrightarrow Y$ such that $f$ is regular over 
the generic point of $Y$ and a very general curve $C$ is 
contracted by $f$ if and only if $D\cdot C=0$. This rational 
map is called the {\em nef reduction} of $D$, and 
$n(X,D):=\dim(Y)$ is called the {\em nef dimension} of $D$.
Our main result is the following

\begin{thm}\label{main} 
Let $(X,B)$ be a log minimal model such that $K+B$ 
has nef dimension at most $3$.
Then the linear system $|m(K+B)|$ is base point free 
for some positive integer $m$.
\end{thm}

As a corollary of Theorem~\ref{main} and the 
Base Point Free Theorem~\cite{KMM}, the $4$-dimensional 
case of the Log Abundance Conjecture is equivalent to 
the following 

\begin{conj}\label{q} Let $(X,B)$ be a log minimal
model with $\dim(X)=4$. If $(K+B)\cdot C> 0$ for very 
general curves $C\subset X$, then $(K+B)^4>0$.
\end{conj}

In fact, Log Abundance and Conjecture~\ref{q} are equivalent 
in dimension $d$ if the Log Minimal Model Program and Log 
Abundance hold for log varieties of dimension at most $d-1$ 
(Theorem~\ref{reduction}). As an application, we give a new
proof of Abundance for smooth minimal surfaces.

Finally, it should be mentioned that the largest class of
singularities for which Log Abundance is expected to hold is 
that of log varieties with {\em log canonical singularities}. 
We expect that Theorem~\ref{k0}, which is the key technical
tool behind the equivalence of Log Abundance and 
Conjecture~\ref{q}, extends to the case when the general 
fibre has log canonical singularities. Therefore the
above mentioned equivalence should be also valid in the case 
of log minimal models with log canonical singularities.


\section{Preliminary}


A {\em variety} is a reduced and irreducible 
separable scheme of finite type, defined over an 
algebraically closed field of characteristic zero. 
A {\em contraction} is a proper morphism $f\colon X
\to Y$ such that $\calO_Y=f_*\calO_X$.


\subsection{Divisors}

Let $X$ a normal variety, and let $L\in \{\Z,\Q,\R\}$. 
An $L$-Weil divisor is an element of $Z^1(X)\otimes_\Z L$.
Two $\R$-Weil divisors $D_1, D_2$ are $L$-{\em linearly 
equivalent}, denoted $D_1\sim_L D_2$, if there exist 
$q_i\in L$ and rational functions $\varphi_i\in k(X)^\times$ 
such that $D_1-D_2=\sum_i q_i(\varphi_i)$. An $\R$-Weil 
divisor $D$ is called
\begin{itemize}
\item[(i)] {\em $L$-Cartier} if $D\sim_L 0$ in 
a neighborhood of each point of $X$.
\item[(ii)] {\em nef} if $D$ is $\R$-Cartier and 
$D\cdot C\ge 0$ for every curve $C\subset X$.
\item[(iii)] {\em ample} if $X$ is projective and 
the numerical class of $D$ belongs to the real cone 
generated by the numerical classes of ample Cartier 
divisors.
\item[(iv)] {\em semi-ample} if there exists a 
contraction $\Phi\colon X\to Y$ and an ample 
$\R$-divisor $H$ on $Y$ such that $D\sim_\R \Phi^*H$. 
If $D$ is rational, this is equivalent to the
linear system $|kD|$ being base point free for some $k$.
\item[(v)] {\em big} if there exists $C>0$ such that
$\dim H^0(X,kD)\ge Ck^{\dim(X)}$ for $k$ sufficiently 
large and divisible. By definition,
$$
H^0(X,kD)=\{a\in k(X)^\times; (a)+kD\ge 0\}\cup \{0\}. 
$$
\end{itemize}

The {\em Iitaka dimension} of $D$ is 
$
\kappa(X,D)=\max_{k\ge 1} \dim \Phi_{|kD|}(X),
$
where $\Phi_{|kD|}\colon X\dashrightarrow {\mathbb P}(|kD|)$ 
is the rational map associated to the linear system 
$|kD|$. If all the linear systems $|kD|$ are empty, 
$\kappa(X,D)=-\infty$. If $D$ is nef, the 
{\em numerical dimension} $\nu(X,D)$ is the largest 
non-negative integer $k$ such that there exists a 
$k$-dimensional cycle $C\subset X$ with 
$D^k\cdot C\ne 0$. 


\subsection{B-divisors}

(V.V. Shokurov~\cite{Logmodels, Plflips}). 
An {\em $L$-b-divisor} $\bD$ of $X$ is a family 
$\{\bD_{X'}\}_{X'}$ of $L$-Weil divisors indexed by 
all birational models of $X$, such that 
$\mu_*(\bD_{X''})=\bD_{X'}$ if $\mu\colon X''\to X'$ 
is a birational contraction. 

Equivalently, $\bD=\sum_E \mult_E(\bD) E$ is a $L$-valued 
function on the set of all geometric valuations of the 
field of rational functions $k(X)$, having finite support 
on some (hence any) birational model of $X$.

\begin{example} (1) Let $\omega$ be a top rational 
differential form of $X$. The associated family of divisors 
$\bK=\{(\omega)_{X'}\}_{X'}$ is called the {\em canonical 
b-divisor} of $X$. 

(2) A rational function $\varphi \in k(X)^\times$ defines a 
b-divisor $\overline{(\varphi)}=\{(\varphi)_{X'}\}_{X'}$.

(3) An $\R$-Cartier divisor $D$ on a birational model $X'$ 
of $X$ defines an $\R$-b-divisor $\overline{D}$ such that
$(\overline{D})_{X''}=\mu^*D$ for every birational contraction 
$\mu\colon X''\to X'$. 
\end{example}

An $\R$-b-divisor $\bD$ is called {\em $L$-b-Cartier} 
if there exists a birational model $X'$ of $X$ such that
$\bD_{X'}$ is $L$-Cartier and $\bD=\overline{\bD_{X'}}$.
In this case, we say that $\bD$ {\em descends to $X'$}.
An $\R$-b-divisor $\bD$ is {\em b-nef}
({\em b-semi-ample}, {\em b-big}, {\em b-nef and good}) 
if there exists a birational contraction $X'\to X$ such 
that $\bD=\overline{\bD_{X'}}$ and $\bD_{X'}$ is 
nef (semi-ample, big, nef and good).


\subsection{Log pairs}

A {\em log pair} $(X,B)$ is a normal variety $X$ 
endowed with a $\Q$-Weil divisor $B$ such that $K+B$ 
is $\Q$-Cartier. A {\em log variety} is a log pair 
$(X,B)$ such that $B$ is effective. The 
{\em discrepancy $\Q$-b-divisor} of a log pair $(X,B)$ 
is the $\Q$-b-divisor of $X$ defined by
$$
\bA(X,B)=\bK-\overline{K+B}.
$$ 
More precisely, fix a top rational differential form 
$\omega\in \wedge^{\dim(X)}\Omega^1_{k(X)/k}$ with 
$K=(\omega)_X$. For a birational contraction 
$\mu\colon Y\to X$, the Weil divisor $(\omega)_Y$ 
is a canonical divisor of $Y$. Then $\bA(X,B)_Y$ is the 
unique $\Q$-Weil divisor on $Y$ such that the following 
adjunction formula holds:
$$
\mu^*((\omega)_X+B)=(\omega)_Y-\bA(X,B)_Y.
$$
It is easy to see that $\bA(X,B)$ is independent of 
the choice of $\omega$ and in fact it is independent 
of the choice of the canonical divisor $K$ in its linear 
equivalence class.

A log pair $(X,B)$ is said to have at most 
{\em Kawamata log terminal singularities} 
if $\mult_E(\bA(X,B))> -1$ for every geometric 
valuation $E$. 


\section{Nef reduction}


The existence of the nef reduction map is 
originally due to Tsuji~\cite{T00}. An algebraic 
proof of the sharper statement below is due to Bauer, 
Campana, Eckl, Kebekus, Peternell, Rams, Szemberg 
and Wotzlaw~\cite{nr}.

\begin{thm}\cite{T00,nr} 
Let $D$ be a nef $\R$-Cartier divisor on a normal
projective variety $X$. Then there exists a rational 
map $f \colon X \dashrightarrow Y$ to a normal projective 
variety $Y$, satisfying the following properties:
\begin{itemize}
\item[(i)]
$f$ is a dominant rational map with connected fibers,
which is a morphism over the general point of $Y$.
\item[(ii)] There exists a countable intersection $U$ of 
Zariski open dense subsets of $X$ such that for every 
curve $C$ with $C\cap U\ne \emptyset$, $f(C)$ is a point
if and only if $D\cdot C=0$.
\end{itemize}
Moreover, $D|_W\equiv 0$ for general fibers $W$
of $f$.
\end{thm}

 The rational map $f$ is unique, and is called the 
{\em nef reduction of } $D$. The dimension of $Y$ is 
called the {\em nef dimension} of $D$, denoted by $n(X,D)$. 
In general, the following inequalities hold~\cite{ab,nr}:
$$
\kappa(X,D)\le \nu(X,D)\le n(X,D)\le \dim(X).
$$

\begin{defn} A nef $\Q$-Cartier divisor $D$ is 
called good if 
$$
\kappa(X,D)=\nu(X,D)=n(X,D).
$$
\end{defn}

\begin{rem} This is equivalent to Kawamata's 
definition~\cite{ab}. If $$\kappa(X,D)=\nu(X,D),$$
there exists a dominant rational map $f\colon 
X\dashrightarrow Y$ and a nef and big $\Q$-divisor 
$H$ on $Y$ such that $\overline{D}\sim_\Q 
f^*(\overline{H})$, by ~\cite{ab}. In particular, 
$n(X,D)$ coincides with the Iitaka and the numerical 
dimension in this case. 
\end{rem}

\begin{rem}\cite{nr} The extremal values 
of the nef dimension are:
\begin{itemize}
\item[(i)] $n(X,D)=0$ if and only if $D$ is 
numerically trivial ($\nu(X,D)=0$).
\item[(ii)] $n(X,D)=\dim(X)$ if and only if 
there exists a countable intersection $U$ of Zariski
open dense subsets of $X$ such that $D\cdot C>0$ for
every curve $C$ with $C\cap U\ne \emptyset$.
\end{itemize}
\end{rem}


\section{Fujita decomposition}


\begin{defn}\cite{Fuj86} 
An $\R$-Cartier divisor $D$ on a normal proper 
variety $X$ has a {\em Fujita decomposition} if 
there exists a b-nef $\R$-b-divisor $\bP$ of $X$
with the following properties:
\begin{itemize}
\item[(i)] $\bP\le \overline{D}$.
\item[(ii)] $\bP=\sup\{\bH; \bH \mbox{ b-nef 
$\R$-b-divisor}, \bH \le \overline{D}\}$.
\end{itemize}
The $\R$-b-divisor $\bP=\bP(D)$ is unique if 
it exists, and is called the {\em semi-positive 
part} of $D$. The $\R$-b-divisor 
$\bE=\overline{D}-\bP$ is called the 
{\em negative part} of $D$, and
$
\overline{D}=\bP+\bE
$
is called the {\em Fujita decomposition} of
$D$.
\end{defn}

\begin{rem}
Allowing divisors with real coefficients is
necessary: there exist Cartier divisors
(in dimension at least $3$) which have a Fujita
decomposition with irrational semi-positive 
part~\cite{Cut}.
\end{rem}

Clearly, a nef $\R$-Cartier divisor $D$ has a 
Fujita decomposition, with semi-positive 
part $\overline{D}$. More examples can be 
constructed using the following property:

\begin{prop}\cite{Fuj86}\label{in} 
Let $f\colon X\to Y$ be a proper contraction, let
$D$ be an $\R$-Cartier divisor on $Y$ and let $E$ 
be an effective $\R$-Cartier divisor on $X$ such that $E$ 
is vertical and supports no fibers over codimension one 
points of $Y$. 

Then $D$ has a Fujita decomposition if and only if 
$f^*D+E$ has a Fujita decomposition, and moreover, 
$\bP(f^*D+E)=f^*(\bP(D))$.
\end{prop}

\begin{lem} Assume LMMP and Log Abundance. Let 
$(X,B)$ be a log variety with log canonical 
singularities. Then $K+B$ has a Fujita decomposition 
if and only if $\kappa(X,K+B)\ge 0$, and the 
semi-positive part is semi-ample. Moreover, 
$$
\bP(K+B)=\overline{K_Y+B_Y},
$$
for a log minimal model $(Y,B_Y)$.
\end{lem}

\begin{lem}\label{vn} \cite{ab,Fuj86}  
Let $f\colon X\to Y$ be a contraction of normal 
proper varieties, and let $D$ be a nef $\R$-divisor
on $X$ which is vertical on $Y$. Then there exists a
b-nef $\R$-b-divisor $\bD$ of $Y$ such that
$\overline{D}=f^*\bD$.
\end{lem}


\section{Lc-trivial fibrations}


This section is a brief introduction to {\em lc-trivial 
fibrations}. We refer the reader to~\cite{bp,modbdiv} for 
more details. As we shall see in a while, these type of 
fibrations appear naturally in inductive arguments in the 
Log Minimal Model Program.

\begin{defn}
An {\em lc-trivial fibration} $f\colon (X,B)\to Y$ 
consists of a contraction of normal varieties $f\colon 
X \to Y$ and a log pair $(X,B)$, satisfying the 
following properties:
\begin{itemize}
\item[(1)] $(X,B)$ has Kawamata log terminal 
singularities over the generic point of $Y$.
\item[(2)] $\rank f_*\calO_X(\lceil \bA(X,B)\rceil)=1$.
\item[(3)] There exist a positive integer $r$, a rational 
function $\varphi\in k(X)^\times$ and a $\Q$-Cartier 
divisor $D$ on $Y$ such that 
$$K+B+\frac{1}{r}(\varphi)= f^*D.$$
\end{itemize}
\end{defn}

For an lc-trivial fibration $f\colon (X,B)\to Y$, we
define $B_Y=\sum_{P\subset Y} b_P P$, where
the sum runs after all prime divisors of $Y$, and  
$$
1-b_P=\sup\{t\in \R; ^\exists U\ni \eta_P, (X,B+tf^*(P)) 
\mbox{ lc sing}/U\}.
$$
It is easy to see that $B_Y$ is a well defined $\Q$-Weil 
divisor on $Y$. By (3), there exists a unique $\Q$-Weil 
divisor $M_Y$ such that the following {\em adjunction 
formula} holds:
$$
K+B+\frac{1}{r}(\varphi)=f^*(K_Y+B_Y+M_Y).
$$
The $\Q$-Weil divisors $B_Y$ and $M_Y$ are called the 
{\em discriminant} and {\em moduli part} of the lc-trivial 
fibration $f\colon (X,B)\to Y$. 

Let now $\sigma\colon Y'\to Y$ be a birational contraction 
from a normal variety $Y'$. Let $X'$ be a resolution of the 
main component of $X\times_Y Y'$ which dominates $Y'$.
The induced morphism $\mu\colon X'\to X$ is birational, and 
let $(X',B_{X'})$ be the crepant log structure on $X'$,
i.e. $\mu^*(K+B)=K_{X'}+B_{X'}$:
\[ \xymatrix{
(X,B) \ar[d]_f & (X',B_{X'}) \ar[l]_{\mu} \ar[d]^{f'}\\
Y        & Y' \ar[l]^{\sigma}
} \]
We say that the lc-trivial fibration 
$f'\colon (X',B_{X'})\to Y'$ is induced by base change. 
Let $B_{Y'}$ be the discriminant of $K_{X'}+B_{X'}$ on 
$Y'$. Since the definition of the discriminant 
is divisorial and $\sigma$ is an isomorphism over
codimension one points of $Y$, we have 
$B_Y=\sigma_*(B_{Y'})$. This means that there exists 
a unique $\Q$-b-divisor $\bB$ of $Y$ such that 
$\bB_{Y'}$ is the discriminant on $Y'$ of the induced 
fibre space $f'\colon (X',B_{X'})\to Y'$, for 
every birational model $Y'$ of $Y$. We call $\bB$ the 
{\em discriminant $\Q$-b-divisor} induced by $(X,B)$ on 
the birational class of $Y$.
Accordingly, there exists a unique 
{\em $\Q$-b-divisor} $\bM$ of $Y$ such that
$$
K_{X'}+B_{X'}+\frac{1}{r}(\varphi)=
f^*(K_{Y'}+\bB_{Y'}+\bM_{Y'})
$$
for every lc-trivial fibration 
$f'\colon (X',B_{X'})\to Y'$ induced by base change on 
a birational model $Y'$ of $Y$. We call $\bM$ the 
{\em moduli $\Q$-b-divisor} of the lc-trivial 
fibration $f\colon (X,B)\to Y$.

The positivity properties of the moduli $\Q$-b-divisor
of an lc-trivial fibration play a key role in applications. 
By~\cite{bp}, Theorem 0.2, the moduli $\Q$-b-divisor 
$\bM$ is b-nef. However, it is expected that $\bM$ is in 
fact b-semiample. This is true, for instance, if $Y$ is 
a curve~\cite{bp}, Theorem 0.1. 
Under an extra assumption, we can prove that $\bM$ is 
``almost'' b-semiample (cf. ~\cite{modbdiv}, Theorem 3.3):

\begin{thm}\label{k0}
Let $f\colon (X,B)\to Y$ be an lc-trivial fibration
such that the geometric generic fibre
$X_{\bar{\eta}}=X\times_Y \Spec(\overline{k(Y)})$ 
is projective and $B_{\bar{\eta}}$ is effective.

Then the moduli $\Q$-b-divisor $\bM$ is b-nef
and good.
\end{thm}


\section{Reduction argument}


\begin{thm}\label{reduction} 
Let $(X,B)$ be a projective log variety with Kawamata 
log terminal singularities such that the log canonical 
divisor $K+B$ is nef, of positive nef dimension 
$n=n(X,K+B)$. 

If the Log Minimal Model Program and Log
Adundance hold in dimension $n(X,K+B)$,
then $K+B$ is a semi-ample $\Q$-divisor.
\end{thm}

\begin{proof} Let $\Phi\colon X \dashrightarrow Y$ be the
quasi-fibration associated to the nef divisor $K+B$, and 
let $\Gamma$ be the normalization of the graph of $\Phi$:
\[ \xymatrix{
 & \Gamma \ar[dl]_\mu \ar[dr]^f & \\
  X      & & Y
} \] 
Since $\Phi$ is a quasi-fibration, $\mu$ is 
birational, $f$ is a contraction and $\Exc(\mu)
\subset \Gamma$ is vertical over $Y$. Let
$W$ be a general fiber of $f$. Let $K_\Gamma+B_\Gamma=
\mu^*(K+B)$ and let $B_W=B_\Gamma\vert_W$.

{\em Step 1}: $(W,B_W)$ is a projective log variety 
with Kawamata log terminal singularities, and 
$K_W+B_W\sim_\Q 0$.
Since $\mu$ is an isomorphism in a neighborhood of
$W$, we infer that $B_W$ is effective, $(W,B_W)$ has
Kawamata log terminal singularities and 
$K_W+B_W=\mu^*(K+B)|_W$. 
The definition of $\Phi$ implies that $K_W+B_W$ is 
numerically trivial. From ~\cite{modbdiv}, 
Theorem 0.1, we conclude that $K_W+B_W\sim_\Q 0$. 

{\em Step 2}: There exist a diagram
\[ \xymatrix{
 X & X' \ar[l]_\mu \ar[d]^{f'}  \\
   &  Y'
} \]
satisfying the following properties:
\begin{enumerate}
\item[(a)] $\mu$ is a birational contraction and
$X',Y'$ are nonsingular.
\item[(b)] $f'\colon (X',B_{X'})\to Y'$ is an 
lc-trivial fibration, where 
$K_{X'}+B_{X'}=\mu^*(K+B)$.
\item[(c)] The moduli $\Q$-b-divisor $\bM$ of the
lc-trivial fibration $f'\colon (X',B_{X'})\to Y'$
descends to $Y'$ and there exists a contraction 
$h\colon Y'\to Z$ and a nef and big $\Q$-divisor 
$N$ on $Z$ such that $\bM_{Y'}\sim_\Q h^*N$.
\item[(d)] Let $E$ be any prime divisor on $X'$.
If $E$ is exceptional over $Y'$, then $E$ is 
exceptional over $X$.
\end{enumerate}

From Step 1, there exists a rational function 
$\varphi\in k(X)^\times$ and $r\in \Q$ such that
$\mu^*(K+B)+r(\varphi)$ is an $f$-vertical nef 
$\Q$-divisor. By Lemma~\ref{vn}, there exists a
b-nef $\Q$-b-divisor $\bD$ of $Y$ such that
$\overline{K+B}+r(\varphi)=f^*(\bD)$.

Let $Y''\to Y$ be a resolution of singularities
such that $\bD$ descends to $Y''$, and let $X''$ be 
the normalization of the main component of 
$X\times_Y Y''$. We denote the induced birational 
contraction by 
$\mu\colon X''\to X$, and let $\mu^*(K+B)=K_{X''}+B_{X''}$.
It is clear that $f'\colon (X'',B_{X''})\to Y''$ is
an lc-trivial fibration. By~\cite{bp}, Theorem 0.2,
the corresponding moduli $\Q$-b-divisor $\bM$ is 
$\Q$-b-Cartier. Furthermore, the $f''$-horizontal
part of $B_{X''}$ is effective since $\Exc(X''/X)$
is $f''$-vertical and $B$ is effective. By~\cite{modbdiv},
Theorem 3.3, we infer that the property (c) holds
if we replace $Y'$ by a sufficiently high resolution.
 
 Let $Y'\to Y''$ be a resolution of singularities
such that (c) holds for the induced lc-trivial fibratian 
and such that $f'\colon X'\to Y'$ dominates a flattening 
of $f$, where now $X'$ is a resolution of singularities 
of the main component of $X\times_Y Y'$. It is clear that 
the properties $(a)-(d)$ hold.

{\em Step 3}: 
There exists an effective $\Q$-divisor $\Delta'$ on 
$Y'$ such that $(Y',\Delta')$ is a log variety with
Kawamata log terminal singularities, $K_{Y'}+\Delta'$
has a Fujita decomposition and 
$$
\overline{K+B}\sim_\Q {f'}^*(\bP(K_{Y'}+\Delta')).
$$

Indeed, let $B_{X'}=B_{X'}^+-B_{X'}^-$ be the decomposition
of $B_{X'}$ into its positive and negative part. There 
exists a unique $\Q$-divisor $F$ on $Y'$ such that 
$A:=B_{X'}^-+{f'}^*F$ is effective and supports no fibers
of $f'$, over a big open subset of $Y'$. In particular,
$f'\colon (X',B_{X'}^+-A)\to Y'$ is an lc-trivial fibration
with the same moduli $\Q$-b-divisor $\bM$. Let $\Delta_{Y'}$
be the corresponding discriminant $\Q$-divisor on $Y'$.
Since $A$ does not support fibers over codimension one points
of $Y'$, we infer that $\Delta_{Y'}$ is effective.
We have 
$$
K_{X'}+B_{X'}^+-A+\frac{1}{b}(\varphi)={f'}^*
(K_{Y'}+\Delta_{Y'}+\bM_{Y'}).
$$
It is clear that $(Y',\Delta_{Y'})$ is a log variety 
with Kawamata log terminal singularities. By (c), 
there exists an effective $\Q$-divisor $\Delta'$ on $Y'$ 
such that $(Y',\Delta')$ is a log variety with Kawamata 
log terminal singularities, and 
$\Delta'\sim_\Q B_{Y'}+\bM_{Y'}$. In particular,
$$
K_{X'}+B_{X'}^+-A \sim_\Q {f'}^*(K_{Y'}+\Delta').
$$
Let $A=A^+-A^-$ be the decomposition of $A$ into
its positive and negative part. We have
$$
\mu^*(K+B)+B_{X'}^-+A^-\sim_\Q {f'}^*(K_{Y'}+\Delta')+A^+.
$$
Both $\Q$-divisors $B_{X'}^-,A^-$ are effective and 
$\mu$-exceptional. This is clear for $B_{X'}^-$. As for
$A^-$, this follows from property (d) since $A^-$ is
$f'$-exceptional by construction. Finally, $A^+$ is
effective, and it does not support fibers over codimension
one points of $Y'$.

Since $K+B$ is nef, the left hand side has a Fujita 
decomposition, with semi-positive part $\overline{K+B}$. 
Proposition~\ref{in} applies, hence $K_{Y'}+\Delta'$ has
a Fujita decomposition as well 
and 
$
\overline{K+B}\sim_\Q {f'}^*(\bP(K_{Y'}+\Delta')).
$

{\em Step 4}: 
From the Log Minimal Model Program and Log Abundance 
applied to the log variety $(Y',\Delta)$, the semi-positive 
part of $K_{Y'}+\Delta'$ is b-semi-ample. Therefore 
$\overline{K+B}$ is b-semi-ample, that is $K+B$ is a 
semi-ample $\Q$-divisor.
\end{proof}

\begin{proof}(of Theorem~\ref{main}) If $K+B$ is numerically 
trivial, then Log Abundance holds by~\cite{modbdiv}, Theorem 
0.1.(1). If $K+B$ has positive nef dimension, we may apply 
Theorem~\ref{reduction}, since the Log Minimal Model Program 
and Log Abundance are known to hold up to dimension $3$ 
(\cite{Sh},\cite{LAB}).
\end{proof}

Finally, we show how Theorem~\ref{reduction} may be 
used to give another proof Abundance in dimension two, in 
characteristic zero.

\begin{thm} Let $X$ be a nonsingular projective 
complex surface whose canonical divisor $K$ is nef. 
Then there exists a positive integer $m$ such that 
the linear system $\vert mK\vert$ is base point free.
\end{thm}

\begin{proof} We may assume that $K$ has maximal nef 
dimension. Indeed, otherwise either $K$ is numerically 
trivial hence a torsion divisor by Kawamata~\cite{minmod},
or $K$ has nef dimension one. The one dimensional
case of Log Abundance is easy to check, hence 
Theorem~\ref{reduction} implies that $K$ is semi-ample if
$n(X,K)=1$.

 If $K$ is big, the result follows from the Base
Point Free Theorem~\cite{KMM}. In the following, we will
assume by contradiction that $K$ has maximal nef dimension
but $(K^2)=0$. In particular, $K\cdot D=0$ for every 
$D\in \vert mK\vert$ ($m\ge 1$). Since $K$ has maximal
nef dimension, we infer $\kappa(K)\le 0$.

Since $(K^2)=0$, the Riemann-Roch formula gives
$$
\chi(X,mK)=1-q(X)+p_a(X).
$$
Furthermore, $H^0(X,mK)=0$ for $m\le -1$, since $K$
has positive nef dimension. Therefore $H^2(X,mK)=0$ 
for $m\ge 2$, by duality. In particular,
$$
h^0(X,mK)\ge 1-q(X)+h^0(X,K) \mbox{ for } m\ge 2.
$$
Assume that $q(X)=0$. In particular, $\kappa(K)\ge 
0$, hence $\kappa(K)=0$. Then there exists a positive
integer $b$ such that
$$
\{m\in \Z_{>0}; \vert mK\vert\ne \emptyset\}=b\Z_{>0}.
$$
But $\vert mK\vert\ne \emptyset$ for every $m\ge 2$, from
above, hence $b=1$. Therefore $h^0(X,K)\ge 1$, hence
$h^0(X,2K)\ge 2$, which contradicts $\kappa(K)=0$.

Assume now that $q(X)>0$. The Stein factorization of the
Albanese map $X\to \Alb(X)$ induces a nontrivial fiber 
space $f\colon X\to Y$. 
Assume first that $Y$ is a curve, and let $F$ be the general 
fibre of $f$. Since $K$ has maximal nef dimension, its
restriction $K_F=K\vert_F$ has maximal nef dimension as well.
Therefore $K_F$ is ample. Furthermore, $\kappa(Y)\ge 0$ 
(\cite{ueno}) and Iitaka's Addition Conjecture $C_{2,1}$ 
(see~\cite{bcv} for instance) implies
$$
\kappa(X)\ge \kappa(F)+\kappa(Y)\ge 1,
$$
a contradiction. Finally, assume that $f$ is birational. 
Since $\kappa(X)\le 0$, we infer by~\cite{ueno}, that $f$ is 
the Albanese map of $X$. Since $K_{\Alb(X)}=0$, there exists 
an effective $f$-exceptional divisor $E$ on $X$ such that $K=E$. 
Therefore $K$ does not have maximal nef dimension, a contradiction.
\end{proof}



\begin{thebibliography}{}


\bibitem{bp}
Ambro, F.,
{\em Shokurov's boundary property,}
{preprint math.AG/0210271.}

\bibitem{abnef}
Ambro, F.,
{\em Nef dimension of minimal models,}
{preprint math.AG/0301305.}

\bibitem{modbdiv}
Ambro, F.,
{\em The moduli b-divisor of an lc-trivial fibration,}
{preprint math.AG/0308143.}

\bibitem{bcv}
Barth, W.; Peters, C.; Van de Ven, A.,
{\em Compact complex surfaces.}
{Ergebnisse der Mathematik und ihrer Grenzgebiete (3), 4.} 
{Springer-Verlag, Berlin, 1984.}

\bibitem{nr}
Bauer, T.; Campana, F.; Eckl, T.; Kebekus, S.; 
Peternell, T.; Rams, S.; Szemberg, T.; Wotzlaw, L.,
{\em A reduction map for nef line bundles.}
{Complex geometry (G\"ottingen, 2000), 27--36, 
Springer, Berlin, 2002.}

\bibitem{Cut}
Cutkosky, S. D.,
{\em Zariski decomposition of divisors on algebraic 
varieties},
{Duke Math. J. {\bf 53} (1986), no. 1, 149--156.} 

\bibitem{Fuj86}
Fujita, T.,
{\em Zariski decomposition and canonical rings of 
elliptic threefolds.}
{J. Math. Soc. Japan {\bf 38} (1986), no. 1, 19--37.}

\bibitem{minmod}
Kawamata, Y.,
{\em Minimal models and the Kodaira dimension of
algebraic fiber spaces},
{J. Reine Angew. Math {\bf 363} (1985), 1-46.}

\bibitem{ab}
Kawamata, Y.,
{\em Pluricanonical systems of minimal algebraic 
varieties},
{Invent. Math., {\bf 79} (1985), 567--588.}

\bibitem{KMM}
Kawamata, Y; Matsuda, K; Matsuki, K,
{\em Introduction to the minimal model program}, 
{Algebraic Geometry, Sendai, Advanced Studies
in Pure Math. {\bf 10} (1987), 283--360.}

\bibitem{Ka92}
Kawamata, Y.,
{\em Abundance theorem for minimal threefolds},
{Invent. Math. {\bf 108} (1992), 229--246.}

\bibitem{LAB}
Keel, S.; Matsuki, K.; McKernan, J., 
{\em Log abundance theorem for threefolds.}
{Duke Math. J. {\bf 75} (1994), no. 1, 99--119.}

\bibitem{LAB'}
Keel, S.; Matsuki, K.; McKernan, J., 
{\em Corrections to~\cite{LAB}}
{Duke Math. J.  {\bf 122} (2004), no. 3, 625--630}

\bibitem{Miy1}
Miyaoka, Y.,
{\em The Chern classes and Kodaira dimension of
a minimal variety. In Algebraic Geometry, Sendai.}
{Adv. Stud. Pure Math. {\bf 10}, Kinokuniya, Tokyo,
1987, 449--476.}

\bibitem{Miy2}
Miyaoka, Y.,
{\em Abundance conjecture for minimal threefolds:
$\nu=1$ case},
{Comp. Math. {\bf 68} (1988), 203--220.}

\bibitem{Miy3}
Miyaoka, Y.,
{\em On the Kodaira dimension of minimal threefolds},
{Math. Ann. {\bf 281} (1988), 325--332.}

\bibitem{reid}
Reid, M.,
{\em Minimal models of canonical $3$-folds.} 
{Algebraic varieties and analytic varieties (Tokyo, 1981),  
131--180, Adv. Stud. Pure Math., 1, North-Holland, Amsterdam, 
1983.}

\bibitem{Sh}
Shokurov, V. V.,
{\em $3$-fold log flips.} 
{Russian Acad. Sci. Izv. Math. {\bf 40} (1993), 
no. 1, 95--202.} 

\bibitem{Logmodels} 
Shokurov, V. V.,
{\em $3$-fold log models.} 
{Algebraic geometry, 4. J. Math. Sci. {\bf 81} (1996),  
no. 3, 2667--2699.} 

\bibitem{Plflips}
Shokurov, V. V.,
{\em Prelimiting flips}, 
{In {\em Birational Geometry: Linear systems and finitely 
generated algebras. A collection of papers.}
Iskovskikh, V.A. and Shokurov, V.V. Editors.
Proc. of Steklov Institute {\bf 240}, 2003.}

\bibitem{T00}
Tsuji, H.,
{\em Numerical trivial fibrations,}
{preprint math.AG/0001023.}

\bibitem{ueno}
Ueno, K.,
{\em Classification of algebraic varieties, I},
{Compositio Math. {\bf 27} (1973), 277--342.}

\end{thebibliography}
\end{document}